\documentclass[12pt,leqno]{article}
\usepackage{amsfonts}
\usepackage{amsmath, amsthm, amsfonts, amssymb, color}
\usepackage{mathrsfs}
\usepackage{color}
\theoremstyle{definition}

\newcommand{\scr}[1]{\mathscr #1}
\definecolor{wco}{rgb}{0.5,0.2,0.3}

\numberwithin{equation}{section} \theoremstyle{remark}

\title{{\bf Closability of Quadratic Forms Associated to Invariant Probability Measures of    SPDEs }\footnote{Supported in part by NNSFC(11131003, 11431014). }}
\author{
{\bf Michael R\"ockner$^{b)}$ and   Feng-Yu Wang$^{a), c)}$}\\
\footnotesize{$^{a)}$School of Mathematical Sciences,
Beijing Normal
University, Beijing 100875, China}\\
\footnotesize{$^{b)}$Department of Mathematics, Bielefeld
University, D-33501 Bielefeld, Germany}\\
 \footnotesize{$^{c)}$Department of Mathematics,
Swansea University, Singleton Park, SA2 8PP, United Kingdom}\\ \footnotesize{wangfy@bnu.edu.cn, F.-Y.Wang@swansea.ac.uk, roeckner@math.uni-bielefeld.de}}

\begin{document}
\allowdisplaybreaks
\def\R{\mathbb R}  \def\ff{\frac} \def\ss{\sqrt}
\def\N{\mathbb N} \def\kk{\kappa}
\def\dd{\delta} \def\DD{\Delta} \def\vv{\varepsilon} \def\rr{\rho}
\def\<{\langle} \def\>{\rangle} \def\GG{\Gamma} \def\gg{\gamma}
  \def\nn{\nabla} \def\pp{\partial} \def\E{\mathbb E}\def\EE{\scr E}
\def\d{\text{\rm{d}}} \def\bb{\beta} \def\aa{\alpha} \def\D{\scr D}
  \def\si{\sigma}
\def\beg{\begin} \def\beq{\begin{equation}}  \def\F{\scr F}
\def\Ric{\text{\rm{Ric}}} \def\Hess{\text{\rm{Hess}}}
\def\e{\text{\rm{e}}}   \def\OO{\Omega}  \def\oo{\omega}
 \def\tt{\tilde}
  \def\P{\mathbb P}
\def\C{\scr C}

\def\Z{\mathbb Z} \def\vrr{\varrho} \def\ll{\lambda}
\def\L{\scr L}
   \def\H{\mathbb H}
\def\M{\scr M}\def\Q{\mathbb Q}   \def\LL{\Lambda}
\def\Rank{{\rm Rank}} \def\B{\scr B}
  \def\BB{\mathbb B}
\maketitle

\begin{abstract} By using the integration by parts formula of a Markov operator, the closability of quadratic forms associated to the corresponding  invariant probability measure is proved. The general result is applied to the study of    semilinear SPDEs, infinite-dimensional stochastic Hamiltonian systems, and semilinear SPDEs with delay.
\end{abstract} \noindent
 \vskip 2cm

\section{Introduction}

Let $\BB$ be a separable Banach space and $\mu$ a reference probability measure on $\BB$. For any $k\in \BB$, let $\pp_k$ denote the directional derivative along $k$. According to \cite{AR}, the form
$$\EE_k(f,g):= \mu((\pp_k f)(\pp_k g)):=\int_\BB (\pp_k f)(\pp_k g)\d\mu,\ \ f,g\in C_b^2(\BB),$$ is closable on $L^2(\mu)$ if
$\rr_s:= \ff{\d\mu(sk+\cdot)}{\d \mu}$ exists for any $s$ such that $s\mapsto\rr_s$ is lower semi-continuous $\mu$-a.e.; i.e. for some fixed $\mu$-versions of $\rr_s, s\in \R$,
$$\liminf_{s\to t}\rr_s(x) \ge \rr_t(x),\ \ \mu{\rm-a.e.}\, x,\ t\in\R.$$
In this paper, we aim to investigate the closability of $\EE_k$ for $\mu$ being the invariant probability measure of a (degenerate/delay) semilinear SPDE. Since in this case the above lower semi-continuity  condition is   hard to check, in this paper we make use of the integration by parts formula for the associated Markov semigroup   in the line of \cite{W12} using coupling arguments.

The main motivation to study the closability of $\EE_k$ (respectively of $\pp_k$) on $L^2(\mu)$ is that it leads to a concept of weak differentiablity on $\BB$ with respect to $\mu$ and one can define the corresponding Sobolev space on $\BB$ in $L^p(\mu)$, $p\in [1,\infty)$. In particular, one can analyze the generator of a Markov process (e.g. arising from a solution of an SPDE) on these Sobolev spaces when $\mu$ is its (infinitesimally) invariant measure, see e.g. \cite{5'} for details.

Before considering  specific models of SPDEs, we first introduce a general result on the closability of $\EE_k$ using the integration by parts formula. To this end, we consider a family of $\BB$-valued random variables $\{X^x\}_{x\in\BB}$ measurable in $x$, and let $P(x,\d y)$ be the distribution of $X^x$ for $x\in\B$.
Then we have the following  Markov operator on $\B_b(\BB):$
$$P f(x):= \int_\BB f(y)P(x,\d y)=\E f(X^x),\ \ x\in\BB, f\in \B_b(\BB).$$ A probability measure $\mu$ on $\BB$ is called an invariant measure of $P$ if
$\mu(Pf)=\mu(f)$   for all $f\in \B_b(\BB).$

\beg{prp}\label{P1.1} Assume that   the Markov operator $P$ has  an invariant probability measure $\mu$. Let $k\in \BB$. If  there exists a family of real random variables $\{M_x\}_{x\in\BB}$ measurable in $x$ such that  $M_\cdot\in L^2(\P\times \mu)$, i.e.
\beq\label{1.0} (\P\times\mu)(|M_\cdot|^2):=\int_\BB \E |M_x|^2 \mu(\d x) <\infty;\end{equation} and the integration by parts formula
\beq\label{1.1} P(\pp_k f)(x)= \E \big\{f(X^x)M_x\big\},\ \ f\in C_b^2(\BB),\mu\text{-a.e.}\ x\in\BB\end{equation}
holds, then $(\EE_k, C_b^2(\BB))$ is closable in $L^2(\mu).$\end{prp}

\beg{proof} Since $\mu$ is $P$-invariant, by \eqref{1.0} and \eqref{1.1} we have
$$ \mu(\pp_k f)= \int_\BB P(\pp_k f)(x)\mu(\d x) = (\P\times \mu)\big( f(X^\cdot) M_\cdot\big),\ \ f\in C_b^2(\BB).$$ So,
\beg{align*}\EE_k(f,g)&:= \mu\big((\pp_k f)(\pp_k g)\big) = \mu(\pp_k\{f\pp_kg\}) - \mu(f\pp_k^2 g)\\
 &= (\P\times \mu)\big(\{f\pp_k g\}(X^\cdot) M_\cdot\big)-\mu(f\pp_k^2g), \ \ f,g\in C_b^2(\BB).\end{align*}  It is standard that this implies the closability of the form $(\EE_k, C_b^2(\BB))$ in $L^2(\mu).$  Indeed, for $\{f_n\}_{n\ge 1}\subset C_b^2(\BB)$ with $f_n\to 0$ and $\pp_k f_n\to Z$ in $L^2(\mu)$, it suffices to prove that $Z=0.$ Since $\mu(f_n^2)\to 0$ and $(\P\times \mu)\big(|f_n\pp_k g|^2(X^\cdot)|\big)=
\mu(|f_n\pp_k g|^2)$ as $\mu$ is $P$-invariant,  the above formula yields
\beg{align*} &|\mu(Z g)|= \lim_{n\to\infty} |\mu(g\pp_k f_n)|\\
&= \lim_{n\to\infty} \big|(\P\times \mu)\big(\{f_n\pp_k g\}(X^\cdot)M_\cdot\big)-\mu(f_n\pp_k^2g)\big|\\
&\le \liminf_{n\to\infty} \Big\{\ss{(\P\times \mu)\big(|f_n\pp_k g|^2(X^\cdot)\big)\cdot (\P\times \mu)(|M_\cdot|^2)}+\ss{ \mu(f_n^2)\mu(|\pp_k^2g|^2)}\Big\}\\
&\le  \liminf_{n\to\infty} \Big\{\|\pp_kg\|_\infty \ss{\mu(f_n^2)\cdot (\P\times \mu)(|M_\cdot|^2)} + \|\pp_k^2g\|_\infty \ss{\mu(f_n^2)}\Big\}=0,\  \ g\in C_b^2(\BB).\end{align*} Therefore, $Z=0$. \end{proof}

\paragraph{Remark 1.1.} The integration by parts formula \eqref{1.1}  implies the estimate
\beq\label{DM} |\mu(\pp_k f)|^2\le (\P\times\mu)(|M_\cdot|^2)\mu(f^2).\end{equation} As the main result in \cite{DD} (Theorem 10), this type of estimate, called Fomin derivative estimate of the invariant measure,   was derived as the main result for the following semi-linear SPDE on $\H:=L^2(\scr O)$ for any bounded open domain $\scr O\subset \R^n$ for $1\le n\le 3$:
$$\d X(t)= [\DD X(t)+p(X(t))]\d t + (-\DD)^{-\gg /2}\d W(t),$$ where $\DD$ is the Dirichlet Laplacian on $\scr O$, $p$ is a decreasing polynomial with odd degree,  $\gg\in (\ff n 2-1, 1),$ and $W_(t)$ is the cylindrical Brownian motion on $\H$.  The main point of the study is to apply the Bismut-Elworthy-Li derivative formula and  the following formula for the semigroup $P_t^\aa$ for the  Yoshida approximation of this SPDE (see \cite[Proposition 7]{DD}):
$$P_t^\aa \pp_k f= \pp_k P_t^\aa -\int_0^t P_{t-s}(\pp_{Ak+\pp_k p}P_s^\aa f) \d s.$$
In this paper we will establish the integration by parts formula of type \eqref{1.1} for the associated semigroup which implies the estimate   \eqref{DM}. Our results  apply to a general framework where the operator $(-\DD)^{-\gg /2}$ is replaced by a suitable linear operator $\si$ (see Section 2)  which can be degenerate (see Section 3), and the drift $p(x)$ is replaced by a general map $b$ which may include a time delay (see Section 4). However, the price we have to pay for the generalization is that the drift $b$ should be regular enough.

\section{Semilinear SPDEs}

Let $(\H,\langle\cdot,\cdot\rangle,|\cdot|)$ be a real
separable Hilbert space,  and  $(W(t))_{t\geq0}$   a cylindrical
Wiener process on $\H$ with respect to  a complete probability space
$(\Omega, \scr {F}, \mathbb{P})$ with the natural filtration
$\{\scr {F}_t\}_{t\geq0}$. Let $\scr {L}(\H)$ and $\scr
{L}_{HS}(\H)$ be the spaces of all linear bounded operators and
Hilbert-Schmidt operators on $H$ respectively. Let $\|\cdot\|$
and $\|\cdot\|_{HS}$  denote the operator norm and the Hilbert-Schmidt norm
respectively.

Consider the following semilinear   SPDE
\beq\label{E1}
\d X(t)=\{AX(t)+b(X(t))\}\d t+\sigma\d W(t),\end{equation}
 where
\begin{enumerate}
\item[ {\bf (A1)}] $(A,\scr {D}(A))$ is a negatively definite self-adjoint  linear
operator on $\H$ with compact resolvent.
\item[ {\bf (A2)}] Let $\H^{-2}$ be the completion of $\H$ under the inner product
$$\<x,y\>_{\H^{-2}}:= \<A^{-1}x, A^{-1}y\>.$$ Let
$b: \H \rightarrow \H^{-2}$ be such that
$$\int_0^1|\e^{tA}b(0)|\d t<\infty,\ \  |\e^{tA}(b(x)-b(y))| \le \gg(t) |x-y|,\ \ x, y\in\H, t>0$$ holds for some positive $\gg\in C((0,\infty))$ with $\int_0^1\gg(t)\d t<\infty.$
\item[ {\bf (A3)}] $ \sigma\in \scr {L}(\H)$ with Ker$(\si\si^*)=\{0\}$ and $\int_0^1 \|\e^{tA} \si\|_{HS}^2\d t<\infty. $  \end{enumerate}

According to {\bf (A1)}, the spectrum of $A$ is discrete with negative eigenvalues. Let $0<\ll_0\le \cdots\le \ll_n\cdots$ be all eigenvalues of $-A$ counting the multiplicities, and let
 $\{e_i\}_{i\ge 1}$ be the corresponding unit eigen-basis. Denote $\H_{A,n}={\rm span}\{e_i: 1\le i\le n\}, n\ge 1$. Then $\H_A:=\cup_{n=1}^\infty \H_{A,n}$  is a dense subspace of  $\H$.
In assumption {\bf (A2)} we have used the fact that for any $t>0$, the operator $\e^{tA}$ extends uniquely to a bounded linear operator from $\H^{-2}$ to  $\H$, which is again denoted by $\e^{tA}.$

Due to assumptions {\bf (A1)}, {\bf (A2)} and {\bf (A3)},  by a standard iteration argument we conclude that for any $x\in\H$   the equation (\ref{E1}) has a unique     mild solution $X^x(t)$ such that $X^x(0)=x$ (see \cite{DZ}).
 Let
 $$P_t f(x)= \E f(X^x(t)),\ \ f\in \B_b(\H), x\in \H$$ be the associated Markov semigroup.

 Let
$$\|x\|_\si= \inf\big\{ |y|: \ y\in\H, \ss{\si\si^*} y=x\big\},\ \ x\in\H,$$ where $\inf\emptyset :=\infty$ by convention.
Then  $\|x\|_\si<\infty$ if and only if $x\in \text{Im}(\si).$

\beg{thm}\label{T2.1} Assume that $P_t$ has an invariant probability measure $\mu$ and $\H_A\subset {\rm Im}(\ss{\si\si^*})$.
\beg{enumerate} \item[$(1)$] For any   $k\in \H_A $ such that
\beq\label{**0}\sup_{x\in\H} \|\pp_k b(x)\|_\si:=\sup_{x\in\H} \limsup_{\vv\downarrow 0} \ff{\|b(x+\vv k)-b(x)\|_\si}\vv <\infty,\end{equation}  the form $(\EE_k, C_b^2(\H))$ is closable in $L^2(\mu)$.
\item[$(2)$] If $\si\si^*$ is invertible and $b:\H\to\H$ is Lipschitz continuous, then $(\EE_k, C_b^2(\H))$ is closable in $L^2(\mu)$ for any $k\in \D(A).$ \end{enumerate}  \end{thm}

\beg{proof} Since $\d\tt W_t:= (\si\si^*)^{-1/2} \si \d W_t$ is also a cylindrical Brownian motion and $\si\d W_t= \ss{\si\si^*}\d \tt W_t$, we may and do assume that $\si$ is non-negatively  definite.

 (1) Without loss of generality,  we may and do assume that $k$ is an eigenvector of $A$, i.e. $Ak= \ll k$ for some $\ll\in \R.$
  We first prove the case where $b$ is Fr\'echet differentiable along the direction $k$.   By $Ak=\ll k$ we have
$$k(t):= \int_0^t \e^{sA}k\d s= \ff{\e^{\ll t}-1}\ll k,\ \ t\ge 0,$$ where for $\ll=0$ we set $\ff{\e^{\ll t}-1}\ll=t$. Due to $\|k\|_\si<\infty$ and \eqref{**0}, the proof of \cite[Theorem 5.1(1)]{W12} leads to the integration by parts formula
\beq\label{B0}P_T(\pp_k f)(x)= \E \big\{f(X^x(T)) M_{x,T}\big\},\ \ f\in C_b^1(\H), x\in\H, T>0,\end{equation} where
$$M_{x,T}:= \ff{\ll}{\e^{\ll T}-1}\int_0^T \Big\<\si^{-1} \Big(k -\ff{\e^{\ll t}-1}\ll (\pp_k b)(X^x(t))\Big),\, \d W(t)\Big\>.$$
Since  \eqref{**0} implies
\beq\label{FH}\int_\BB \E|M_{x,T}|^2\,\mu(\d x) \le  \ff{\ll^2}{(\e^{\ll T}-1)^2}\int_0^T\Big\|\si^{-1} \Big(k -\ff{\e^{\ll t}-1}\ll \pp_k b\Big)\Big\|_\infty^2\d t<\infty,\end{equation}
 $(\EE_k, C_b^2(\H))$ is closable in $L^2(\mu)$ according to   Proposition \ref{P1.1}.

 In general, for any $\vv>0$ let
$$b_\vv(x)= \ff 1 {\ss{2\pi\vv}} \int_{\R} b(x+rk) \exp\Big[  -\ff{r^2}{2\vv}\Big] \d r,\ \ x\in \H. $$ Then for any $\vv>0$, $b_\vv$ is Fr\'echet differentiable along $k$ and \eqref{**0} holds uniformly in $\vv$ with $b_\vv$ replacing $b$.    Let $P_t^\vv$ be the semigroup for the solution $X_\vv(t)$ associated to equation \eqref{E1} with $b_\vv$ replacing $b$. By  simple calculations we have:
\beg{enumerate}
\item[(i)] $\lim_{\vv\downarrow 0}\E|X_\vv^x(t)-X^x(t)|^2 =0,    \ t\ge 0, x\in \H.$
\item[(ii)] For any $T>0$, the family
$$M_{\cdot,T}^\vv:=  \ff{\ll}{\e^{\ll T}-1}\int_0^T \Big\<\si^{-1} \Big(k -\ff{\e^{\ll t}-1}\ll (\pp_k b_\vv)(X^\cdot_\vv(t))\Big),\, \d W(t)\Big\>,\ \ \ \vv>0$$ is bounded in $L^2(\P\times\mu)$; i.e. $\sup_{\vv>0} \int_\BB \E |M_{x,T}|^2\,\mu(\d x)<\infty.$
\item[(iii)] $P_T^\vv(\pp_k f)(x)= \E \big(f(X^x_\vv(T) M_{x,T}^\vv\big),\  f\in C_b^1(\H), \vv>0.$\end{enumerate}
So, there exist $M_{\cdot,T}\in L^2(\P\times\mu)$ and  a sequence $\vv_n\downarrow 0$ such that $M_{\cdot,T}^{\vv_n}\to M_{\cdot,T}$ weakly in $L^2(\P\times\mu)$. Thus,  by taking $n\to\infty$ in (iii) and using (i), we prove \eqref{B0} for $\mu$-a.e. $x\in\BB$. Then the   proof of the first assertion  is completed as in the first case.

(2) Since $\si$ is invertible, {\bf (A3)} implies $\aa:= \sum_{i=1}^\infty \ff 1 {\ll_i}<\infty.$ Next, since
 the Lipschitz constant $\|\pp b\|_\infty$ of $b$ is finite,  the integration by parts formula \eqref{B0} also implies explicit Fomin derivative estimates on the invariant probability measure, which were investigated recently in \cite{DD}. Indeed,  it follows from    \eqref{B0} and \eqref{FH} that
\beg{align*}  |\mu(\pp_k f)| &= \inf_{T>0} |\mu(P_T (\pp_k f))|\le \inf_{T>0} \ss{\mu(P_Tf^2)} \bigg(\int_\BB \E |M_{x,T}|^2\,\mu(\d x)\bigg)^{\ff 1 2}\\
&\le |k|\cdot \|f\|_{L^2(\mu)} \inf_{T>0}  \ff{\ll}{\e^{\ll T}-1}\bigg(\int_0^T\Big\|\si^{-1} \Big(I -\ff{\e^{\ll t}-1}\ll \pp b\Big)\Big\|_\infty^2\d t\bigg)^{\ff 1 2},\ \ Ak=\ll k.\end{align*}
By taking $k=e_i, T=\ll_i^{-1}$ and  $\ll=-\ll_i$ in the above estimate,   for any $k\in\D(A)$ we have
\beq\label{*Z}\beg{split} |\mu(\pp_k f)| &\le \sum_{i=1}^\infty |\<k, e_i\> \mu(\pp_{e_i}f)| \le   \bigg(\sum_{i=1}^\infty \ll_i^2 \<k,e_i\>^2 \bigg)^{\ff 1 2}
 \bigg(\sum_{i=1}^\infty \ff 1 {\ll_i^2} \mu(\pp_{e_i} f)^2   \bigg)^{\ff 1 2}\\
 &\le |Ak|  \bigg(\sum_{i=1}^\infty  \ff{\|\si^{-1}\|^2}{\ll_i (\e-1)^2} \Big(1+ \ff{\e-1}{\ll_i}\|\pp b\|_\infty\Big)^2\bigg)^{\ff 1 2}\|f\|_{L^2(\mu)}\\
 &\le C|Ak|\cdot\|f\|_{L^2(\mu)},\end{split}\end{equation}  where
 $C:= \ff{\|\si^{-1}\|\ss{\aa}}{\e-1}\Big(1+ \ff{\e-1}{\ll_1} \|\pp b\|_\infty\Big).$ This implies the closablity of $(\EE_k, C_b^2(\H))$ as explained in the proof of Proposition \ref{P1.1}.
 Indeed, if $\{f_n\}_{n\ge 1}\subset C_b^2(\BB)$ satisfies $f_n\to 0$ and $\pp_k f_n\to Z$ in $L^2(\mu)$, then \eqref{*Z} implies
 \beg{align*} |\mu(gZ)|&= \lim_{n\to\infty} |\mu(g\pp_k f_n)| =\lim_{n\to\infty} |\mu(\pp_k(f_n g)- \mu(f_n\pp_k g)|\\
 &\le C|Ak|\lim_{n\to\infty}  \ss{\mu((f_ng)^2)}=0,\ \ \ g\in C_b^2(\BB),\end{align*} so that $Z=0$. \end{proof}

\

To conclude this section, let us recall a result concerning existence and stability of the invariant probability measure. Let $W_a(t)=\int_0^t \e^{A(t-s)}\si\d W(s), t\ge 0.$  Assume that $b$ is Lipschitz continuous and $\int_0^\infty\|\e^{tA}\si\|_{HS}^2\d t<\infty$. We have
$$\sup_{t\ge 0} \E \big(\|W_A(t)\|^2  + |b(W_A(t))|^2\big)<\infty.$$ Therefore, by \cite[Theorem 2.3]{DZ2}, if there exist $c_1>0,c_2\in\R$ with $c_1+c_2>0$ such that
$$\<A(x-y), x-y\>\le -c_1 |x-y|^2,\ \<b(x)-b(y),x-y\>\le -c_2|x-y|^2,\ \ x,y\in\H,$$ then $P_t$ has a unique invariant probability measure such that $\lim_{t\to\infty} P_t f=\mu(f)$ holds for $f\in C_b(\H).$

\section{Stochastic Hamiltonian systems on Hilbert spaces}

Let $\tt\H$ and $\H$ be two separable Hilbert spaces. Consider the following stochastic differential equation for $Z(t):= (X(t), Y(t))$ on $\tt\H\times\H$:
 \beq\label{3.1} \beg{cases} \d X(t)=   BY(t) \d t,\\
\d Y(t)= \{AY(t)+ b(t,X(t),Y(t))\}\d t +\si \d W(t),\end{cases}\end{equation} where   $B\in\scr L(\H\to \tt\H)$,  $(A,\D(A))$ satisfies {\bf (A1)}, $\si$ satisfies {\bf (A3)},   $W(t)$ is the cylindrical Brownian motion on $\H$, and $b: [0,\infty)\times \tt\H\times\H\to \H^{-2}$  satisfies:  for any $T>0$ there exists   $\gg\in C((0,T])$ with $\int_0^T\gg(t)\d t<\infty$ such that
\beq\label{BB} \beg{split} &\sup_{s\in [0,T]}\int_0^T|\e^{tA} b(s,0)|\d t<1,\\
&\sup_{s\in [0,T]}|\e^{tA}(b(s,z)-b(s,z'))|\le\gg(t)|z-z'|,\ \ t\in [0,T], z,z'\in\tt\H\times\H.\end{split}\end{equation}
Obviously, for any initial data $z:=(x,y)\in \H$, the equation has a  unique mild solution $Z^z(t).$ Let $P_t$ be the associated Markov semigroup.

When $\tt\H$ and $\H$ are finite-dimensional, the integration by parts formula of $P_t$ has been established in \cite[Theorem 3.1]{W12}. Here, we extend this result to the present infinite-dimensional setting.

\beg{prp}\label{P3.1} Assume that $BB^*\in \scr L(\tt\H)$ with ${\rm Ker}(BB^*)=\{0\}.$ Let $T>0$ and $k:=(k_1,k_2)\in {\rm Im}(BB^*)\times\H$ be such that
\beq\label{EG}A k_2=\theta_2 k_2,\ \  A B^*(BB^*)^{-1} k_1= \theta_1B^*(BB^*)^{-1} k_1\end{equation} for some constants $\theta_1,\theta_2\in \R$. For any $\phi,\psi\in C^1([0,T])$ such that
\beq\label{PH} \phi(0)=\phi(T)=\psi(0)=\psi(T)-1=\int_0^T\e^{\theta_2 t}\psi(t)\d t=0,\ \ \int_0^T \phi(t)\e^{\theta_1 t}\d t=\e^{\theta_1 T}, \end{equation} let
\beg{equation*}\beg{split}
&h(t)=B^*(BB^*)^{-1} k_1 \int_0^t \phi'(s) \e^{\theta_1(s-T)}\d s+ k_2 \int_0^t\psi'(s) \e^{\theta_2(s-T)}\d s,\\
&\tt h(t)= \phi(t) \e^{\theta_1(t-T)} B^*(BB^*)^{-1} k_1 +\psi(t) \e^{\theta_2(t-T)} k_2,\\
&\Theta(t) =\bigg(\int_0^t B\tt h(s)\d s,\ \tt h(t)\bigg),\ \  \ t\in [0,T].\end{split}\end{equation*}
If for any $t\in [0,T]$, $b(s,\cdot)$ is Fr\'echet differentiable along $\Theta(t)$  such that
\beq\label{*A} \int_0^T \sup_{z\in \tt\H\times\H}\big\|h'(t)-(\pp_{\Theta(t)}b(t,\cdot))(z)\big\|_\si^2 \d t<\infty,\end{equation} then for any $f\in C_b^1(\tt\H\times\H),$
$$P_T(\pp_k f)= \E\bigg\{f(Z(T))\int_0^T\Big\<(\si\si^*)^{-1/2}\big\{h'(t)-(\pp_{\Theta(t)}b(t,\cdot))(Z(t))\big\},\ \d W(t)\Big\>\bigg\}.$$
\end{prp}

\beg{proof} As explained in the proof of Theorem \ref{T2.1}, we simply assume that $\si=\ss{\si\si^*}$.
Let $(X^0(t),Y^0(t))=(X(t),Y(t))$ solve (\ref{3.1}) with initial data $(x,y)$, and for  $\vv\in (0,1]$ let    $(X^\vv(t), Y^\vv(t))$ solve the equation
\beq\label{CCC0}\beg{cases} \d X^\vv(t) = B Y^\vv(t) \d t,\ \ X^\vv(0) =x,\\
\d Y^\vv(t) = \si \d W(t) +\big\{b(t,X(t),Y(t))+AY^\vv(t)+\vv h'(t)\big\}\d t,\ \ Y^\vv(0)=y.\end{cases}\end{equation} Then it is easy to see from \eqref{EG} and \eqref{PH} that
\beg{equation*}  \beg{split}  &Y^\vv(t)-Y(t) =\vv \int_0^t \e^{(t-s)A}h'(s)\d s\\
&=\vv B^*(BB^*)^{-1} k_1 \int_0^t \phi'(s)\e^{\theta_1(s-T)}\e^{\theta_1(t-s)}\d s +\vv k_2 \int_0^t \psi'(s)\e^{\theta_2(s-T)}\e^{\theta_2(t-s)}\d s\\
&=\vv\big(\phi(t)\e^{\theta_1(t-T)} B^*(BB^*)^{-1} k_1 +\psi(t) \e^{\theta_2(t-T)} k_2\big)=\vv\tt h(t),\end{split}\end{equation*}
and hence,
\beg{equation*}  \beg{split} &X^\vv(t)-X(t)=\vv\int_0^t B\tt h(s)\d s\\
& = \vv \bigg(k_1\int_0^t\phi(r)\e^{\theta_1(r-T)}\d r + (Bk_2) \int_0^t \psi(r) \e^{\theta_2(r-T)} \d r\bigg).\end{split} \end{equation*}So,
\beq\label{X3} X^\vv(t)-X(t)=\vv \Theta(t),\ \ \ t\in [0,T],\end{equation} and in particular
\beq\label{CP} (X^\vv(T),Y^\vv(T))= (X(T),Y(T))+\vv k\end{equation} due to  \eqref{PH}.  Next, \beq\label{XI} \xi_\vv(s)= \vv h'(s)+ b(s, X(s),Y(s))-b(s, X^\vv(s), Y^\vv(s))\end{equation} and
$$R_\vv=\exp\bigg[- \int_0^T \big\<\si^{-1}\xi_\vv(s),\d W(s)\big\>-\ff {1} 2 \int_0^T |\si^{-1}\xi_\vv(s)|^2\d s\bigg].$$  We reformulate (\ref{CCC0}) as
\beq\label{CCC}\beg{cases} \d X^\vv(t) = B Y^\vv(t) \d t,\ \ X^\vv(0) =x,\\
\d Y^\vv(t) = \si  \d W^\vv(t) + \{b(t, X^\vv(t),Y^\vv(t))+AY^\vv(t)\} \d t,\ \ Y^\vv(0)=y,\end{cases}\end{equation} where  by \eqref{*A} and \eqref{X3},
$$W^\vv(t):= W(t)+\int_0^t \si^{-1}\xi_\vv(s)  \d s,\ \ t\in [0,T]$$
is a cylindrical Brownian motion under the weighted probability measure $\Q_\vv:=R_\vv\P$.
Since    $|\xi_\vv|$ is uniformly bounded on $[0,T]$,  by the dominated convergence theorem and \eqref{X3},    for any $f\in C_b^1(\tt\H\times\H)$
we obtain
\beg{equation*}\beg{split} P_T(\pp_kf) &= \lim_{\vv\to 0} \E \ff{f((X(T),Y(T))+\vv k) -f((X(t),Y(t)))}\vv \\
&= \lim_{\vv\to 0} \E \ff{f((X^\vv(T),Y^\vv(T))) -R_\vv f((X^\vv(T),Y^\vv(T)))}\vv\\
&= \E\bigg(f(Z(T)) \lim_{\vv\to 0} \ff{1-R_\vv}\vv\Big)\\
 &= \E\bigg(f(Z(T))\int_0^T\Big\<\si^{-1}\big\{h'(t)-(\pp_{\Theta(t)}b)(Z(t))\big\},\ \d W(t)\Big\>\bigg).\end{split}\end{equation*}
\end{proof}

To apply  this result, we present here a specific choice  of ($\phi,\psi$) such that    \eqref{PH} holds:
$$\phi(t)= \ff{\e^{\theta_1 T}t(T-t)}{\int_0^T s(T-s)\e^{\theta_1 s}\d s},\ \ \psi(t)= \ff{\e^{\theta_2(T-t)}}T\Big(\ff{3t^2}T-2t\Big),\ \ \ \ t\in [0,T].$$

\beg{thm}\label{T3.1} Let $\tt\H=\H=\H$ and ${\rm Ker} (B)=\{0\}.$  Let $b(t,\cdot)=b$ do not dependent on $t$ such that $P_t$ has an invariant probability measure $\mu$.
 If
 \beq\label{*A2}   \sup_{(x,y)\in\H\times\H}\lim_{r\downarrow 0} \ff{\|b(x+r B^{-1}\tt k,y+r k)-b(x,y)\|_\si}r<\infty,\ \ (\tt k,k)\in   (B\H_A)\times\H_A,\end{equation} Then  for any
 $(k_1,k_2)\in (B\H_A)\times\H_A$,   the form $(\EE_k, C_b^2(\H\times\H))$ is closable in $L^2(\mu).$
\end{thm}

\beg{proof} It suffices to prove for $k=(k_1,k_2)$ such that $B^{-1}k_1$ and $k_2$ are eigenvectors of $A$, i.e. $AB^{-1}k_1=\theta_1 B^{-1}k_1$ and $ Ak_2=\theta_2k_2$
hold for some $\theta_1,\theta_2\in\R$. As explained above there exists $T>0$ such that \eqref{PH} holds for some $\phi,\psi\in C^\infty([0,T])$. Moreover, as explained  in the proof of Theorem \ref{T2.1}, by taking
$$b_\vv(s,x,y)= \ff 1 {\ss{2\pi\vv}}\int_\R b\big((x,y)+r \Theta(s)\big)\exp\Big[-\ff{r^2}{2\vv}\Big]\d r,\ \ s\in [0,T],  (x,y)\in\H\times\H $$ for $\vv>0$, such that \eqref{*A2} holds uniformly in $\vv>0$ and $s\in [0,T]$ with $b_\vv(s,\cdot)$ replacing $b$,
 we may and do assume that $b(s,\cdot)$ is Fr\'echet differentiable along $\Theta(s)$. Then  the integration by parts formula in Proposition \ref{P3.1} holds,
 and due to \eqref{*A2} we have
 $$M_{\cdot, T}:=\int_0^T\Big\<(\si\si^*)^{-1/2}\big\{h'(t)-(\pp_{\Theta(t)}b(t,\cdot))(Z(t))\big\},\ \d W(t)\Big\>\in L^2(\P\times\mu).$$ Therefore,  by Proposition  \ref{P1.1},  the form  $(\EE_k, C_b^2(\H\times\H))$ is closable on $L^2(\mu).$
\end{proof}

Below are   typical examples of  the stochastic Hamiltonian system with invariant probability measure such that Theorem \ref{T3.1} applies.

\paragraph{Example 3.1.} Let $\tt\H=\H=\H$.

(1) Let $\H=\R^d$ for some $d\ge 1$. When $\si=B=I$, $A\le-\ll I$ for some $\ll>0$ is a negatively definite $d\times d$-matrix,   and $b(x,y)=A^{-1}\nn V(x)$ for some $V\in C^2(\R^d)$ such that
$\int_{\R^d} \e^{-V(x)}\d x<\infty.$ Then the unique invariant probability measure of $P_t$ is
$$\mu(\d x,\d y)= C \e^{-V(x)+\ff \ll 2 \<Ay,y\>}\d x\d y,$$ where $C>0$ is the normalization. See \cite{Baudoin, GR, V} for the study of hypercoercivity of the associated semigroup $P_t$ with respect to $\mu$, as well as \cite{W15} for the stronger property of  hypercontractivity.

(2) In the infinite-dimensional setting,  let $\si=B=I$  and  $A$ be negatively definite such that $A^{-1}$ is of trace class. Take  $b(x,y)=A^{-1}Qx$ for some positively definite self-adjoint operator $\Q$ on $\H$ such that $Q^{-1}$ is of trace class and
$$\int_0^1\|\e^{t A}A^{-1}Q\|\d t<1.$$  Then it is easy to see that
 $$\mu(\d x,\d y)= N_{Q^{-1}}(\d x) N_{-A^{-1}}(\d y)$$ is an invariant probability measure.

 (3) More generally, let $\si=B=I$ and
 $$b(x,y)= \tt b(x):= A^{-1}\nn V(x),\ \ (x,y)\in \H\times \H_A$$   for some Fr\'echet differentiable $V: \H_A\to \R$ such that \eqref{*A2} holds.   For any $n\ge 1$, let
 $$V_n(r)= V\circ\varphi_n(r),  \ \varphi_n(r)=\sum_{i=1}^n   r_i e_i,\ \ r=(r_1,\cdots, r_n)\in\R^n.$$  If $\int_{\R^n}\e^{-V_n(r)}\d r<\infty$ and when $n\to\infty$   the probability measure
 $$\nu_n(D):= \ff 1{\int_{\R^n}\e^{-V_n(r)}\d r}  \int_{\varphi_n^{-1}(D)} \e^{-V_n(r)}\d r,\ \ D\in \B(\H)$$  converges weakly to some probability measure $\nu$,   then
 $\mu:= \nu\times  N_{-A^{-1}}$ is an invariant probability measure of $P_t$. This can be confirmed by (1) and a finite-dimensional approximation argument. Indeed, let $\pi_n:\H\to \H_{A,n}$ be the orthogonal projection, and let $ A_n= \pi_n A,   W_n= \pi_n W$ and
  $b_n(x,y)= \pi_n \nn V(x)$. Let $X_n(t)$ solve the finite-dimensional equation
 \beg{equation*}\beg{cases} \d X_n(t)=   Y_n(t) \d t,\\
\d Y_n(t)= \{A_nY_n(t)+ b_n(X_n(t))\}\d t +\d W_n(t)\end{cases}\end{equation*} with $(X_n(0), Y_n(0))= (\pi_n X(0), \pi_n Y(0)).$ Then the proof of \cite[Theorem 2.1]{WZ} yields that for every $t\ge 0$,
$$\lim_{n\to \infty} \E \big(|X_n(t)-X(t)|^2+|Y_n(t)-Y(t)|^2\big)=0$$ uniformly in the initial data $(X(0),Y(0))\in\H\times\H$. Thus, letting $P_t^{(n)}$ be the semigroup for $(X_n(t), Y_n(t))$, we have
$$\lim_{n\to\infty} \sup_{(x,y)\in\H\times\H} |P_t^{(n)} f(\pi_n x,\pi_n y)-P_tf(x,y)|=0,\ \ f\in C_b^1(\H\times\H).$$ Combining this with the assertion in (1) and noting that $\nu_n\times (N_{-A^{-1}}\circ\pi_n^{-1})\to\mu$ weakly as $n\to\infty$, we conclude that $\mu$ is an invariant probability measure of $P_t$.

\section{Semilinear SPDEs with delay}

For fixed $\tau>0$, let $\C_\tau= C([-\tau,0];\H)$ be equipped with the uniform norm $\|\eta\|_\infty:=\sup_{\theta\in [-\tau,0]}|\eta(\theta)|.$
For any $\xi\in C([-\tau, \infty);\H)$, we define $\xi_\cdot\in C([0,\infty);\C_\tau)$ by letting
$$\xi_t(\theta)= \xi(t+\theta),\ \   \theta\in [-\tau,0], t\ge 0.$$
Consider the following stochastic   differential   equation with delay:
\beq\label{4.1}\d X(t)= \big\{A X(t)+b(X_t)\big\}\d t+\si\d W(t),\ \ X_0\in\C_\tau, \end{equation}
where $(A,\D(A))$ satisfies {\bf (A1)}, $\si$ satisfies {\bf (A3)}, and $b: \C_\tau\to \H$ satisfies: for any $T>0$ there exists $\gg\in  C((0,T]) $ with $\int_0^T\gg(t)\d t<\infty$ such that
\beq\label{*4} \int_0^T \sup_{s\in [0,T]} |\e^{tA} b(s,0)|^2\d t<\infty, \ \ |\e^{tA} (b(s,\xi)-b(s,\eta))|^2\le \gg(t)\|\xi-\eta\|_\infty^2,\ \ t,s\in [0,T].\end{equation}
 Then     for any initial datum
$\xi\in \C_\tau$, the equation has a unique mild solution $X^\xi(t)$ with $X_0=\xi$. Let $P_t$ be the Markov semigroup for the segment solution $X_t.$

Let
$$  \C_\tau^1 =\bigg\{\eta\in \C_\tau:  \eta(\theta)\in \D(A) \ \text{for}\ \theta\in [-\tau,0],
\int_{-\tau}^0 \big(|A\eta(\theta)|^2+ |\eta'(\theta)|^2\big)\d\theta <\infty\bigg\}.$$
The following result is an extension of \cite[Theorem 4.1(1)]{W12} to the infinite-dimensional setting.

\beg{prp}\label{P4.1}
For any $\eta\in \C_\tau^1$ and $T>\tau$, let $$\GG(t):= \beg{cases} \ff 1 {T-\tau} \e^{(s+\tau-T)A}\eta(-\tau), &\text{if}\ s\in [0,T-\tau],\\
\eta'(s-T)-A\eta(s-T), &\text{if}\ s\in (T-\tau,T],\end{cases}$$
and $$\Theta(t):= \int_0^{t\lor 0} \GG(s)\d s,\ \ t\in [-\tau,T].$$ If  $b(t,\cdot)$ is Fr\'echet differentiable along $\Theta_t$ for $t\in [0,T]$ such that
\beq\label{SS}\sup_{\xi\in\C_\tau} \int_0^T\big\|\GG(t)-(\nn_{\Theta_t}b(T,\cdot))(\xi)\big\|_\si^2 \d t<\infty,\end{equation} then
\beq\label{PP}P_T(\pp_\eta f)=\E\bigg(f(X_T)\int_0^T\Big\<(\si\si^*)^{-1/2}\big(\GG(t)-(\nn_{\Theta_t}b(t,\cdot))(X_t)\big), \d W(t)\Big\>\bigg),\ \ f\in C_b^1(\C_\tau).\end{equation}
 \end{prp}

\beg{proof} Simply let $\si=\ss{\si\si^*}$ as in the proof of Theorem \ref{T2.1}. For any $\vv\in (0,1)$, let $X^\vv(t)$ solve the equation
\beq\label{E'}\d X^\vv(t)= \{AX^\vv(t)+b(t,X_t)+\vv\GG(t)\}\d t+\si \d W(t),\ \ X^\vv_0=X_0.\end{equation} We have
\beg{equation}\label{AP} \beg{split} &X^\vv(t)-X(t)= \vv\int_0^{t^+} \e^{(t-s)A} \GG(s)\d s \\
&= \ff{\vv t^+}{T-\tau} \e^{(\tau-T)A} \eta(-\tau) 1_{[-\tau, T-\tau)}(t) + \vv\eta(t-T)1_{[T-\tau,T]}(t),\ \ t\in [-\tau,T].\end{split}\end{equation} In particular, we have $X_T^\vv-X_T=\vv\eta.$ To formulate $P_T$ using $X_T^\vv$,   rewrite \eqref{E'} by
$$\d X^\vv(t)= \{AX^\vv(t)+b(t,X_t^\vv)\}\d t+\si \d W_\vv(t),\ \ X^\vv_0=X_0,$$where
$$W_\vv(t):= W(t)+\int_0^t \xi_\vv(s)\d s,\ \
\xi_\vv(s):=  b(s,X_s)-b(s,X_s^\vv) +\vv \GG(s).$$ By \eqref{SS} and the Girsanov theorem, we see that $\{W_\vv(t) \}_{t\in [0,T]}$ is a cylindrical Brownian motion on $\H$ under the probability measure $\d\Q_\vv:= R_\vv\d\P$, where
$$R_\vv:=\exp\bigg[\int_0^T\Big\<\si^{-1}\big(b(t,X_t^\vv)-b(t,X_t)-\vv\GG(t)\big),\ \d W(t)\Big\>\bigg].$$ Then
$$\E(f(X_T))= P_Tf= \E (R_\vv f(X_T^\vv)).$$ Combining this with $X_T^\vv= X_T+\vv\eta$ and using \eqref{AP}, we arrive at
\beg{equation*}\beg{split} &P_T(\pp_\eta f)= \lim_{\vv\downarrow  0}  \ff 1 \vv \E\{f(X_T+\vv\eta)-f(X_T)\}=\lim_{\vv\downarrow  0}  \ff 1 \vv \E\{f(X_T^\vv)-R_\vv f(X_T^\vv)\}\\
&= \E\Big(f(X_T)\lim_{\vv\downarrow  0}\ff{1-R_\vv}\vv\Big)= \E\bigg\{f(X_T)\int_0^T \Big\<\si^{-1}\big(\GG(t)-(\nn_{\Theta_t}b(t,\cdot))(X_t)\big), \d W(t)\Big\>\bigg\}.\end{split}\end{equation*}
\end{proof}

\beg{thm}\label{T4.1}  Let $b(t,\cdot)=b$ be independent of $t$ such that $P_t$ has an invariant probability measure $\mu$.  If ${\rm Im}(\si)\supset\H_A$ and
 \beq\label{CD} \sup_{\xi\in\C_\tau}\limsup_{\vv\downarrow 0}  \ff{\|b(\xi+\vv\eta)-b(\xi)\|_\si}{\vv} <\infty,\ \ \eta\in \C_\tau^1\cap \Big(\cup_{n\ge 1}C([-\tau,0];\H_{A,n})\Big),\end{equation} then for any $\eta\in \C_\tau^1\cap \big(\cup_{n\ge 1}C([-\tau,0];\H_{A,n})\big)$, which is dense in $\C_\tau$, the form
$$\EE_\eta(f,g):= \int_{\C_\tau} (\pp_\eta f)(\pp_\eta g)\d\mu,\ \ f,g\in C_b^2(\C_\tau)$$ is closable in $L^2(\mu)$. \end{thm}

\beg{proof} For any $\vv\in (0,1)$ let
$$b_\vv(t,\xi)= \ff 1 {\ss{2\pi\vv}}\int_\R b(\xi+r\Theta_t)\exp\Big[-\ff{r^2}{2\vv}\Big] \d r,\ \ \xi\in \C_\tau.$$ Then $b_\vv(t,\cdot)$ is F\'echet differentiable along $\Theta_t$ and \eqref{CD} holds uniformly in $\vv$ with $b_\vv(t,\cdot)$ replacing  $b$. Moreover, $\eta\in \C_\tau^1\cap \big(\cup_{n\ge 1}C([-\tau,0];\H_n)\big)$ implies that $\Theta_t\in \C_\tau^1\cap \big(\cup_{n\ge 1}C([-\tau,0];\H_n)\big)$ and \eqref{CD} holds uniformly in $t\in [0,T]$ and $\vv\in (0,1)$  with $\Theta_t$ and $b_\vv(t,\cdot)$ replacing $\eta$ and $b$ respectively.  Combining this with    Im$(\si)\supset \H_A$, we conclude that \eqref{SS} holds uniformly in $\vv$  with $b_\vv$ replacing  $b$. Therefore,
  as explained in the proof of Theorem \ref{T2.1}, we may assume that $b$ is Fr\'echet differentiable along $\Theta_t, t\in [0,T],$  and by   Proposition \ref{P4.1} the integration by parts formula \eqref{PP} holds. Moreover, \eqref{CD} implies
 $$M_{\cdot,T}:=  \int_0^T\Big\<(\si\si^*)^{-1/2}\big(\GG(t)-(\nn_{\Theta_t}b(t,\cdot))(X_t)\big), \d W(t)\Big\>\in L^2(\P\times\mu).$$
   Then the proof is finished by Proposition \ref{P1.1}. \end{proof}

Finally,  we introduce the following example  to illustrate Theorem \ref{T4.1}.

\paragraph{Example 4.1.}  Let $b(\xi)= F(\xi(-\tau)), \xi\in\C_\tau,$ for some $F\in C_b^1(\H)$. If $\si$ is Hilbert-Schmidt and
$$\<x, Ax +F(y)-F(y')\>\le -\ll_1|x|^2+\ll_2|y-y'|^2,\ \ x,y\in\H,$$    for some constants $\ll_1>\ll_2\ge 0$, then according to \cite[Theorem 4.9]{BTY} $P_t$ has a unique invariant probability measure $\mu$.
If moreover Im$(\si)\supset \H_A$ and for any $y\in\H_A$ there exists a constant
$$\limsup_{\vv\downarrow 0} \sup_{x\in\H}\ff{\|F(x+\vv y)-F(x)\|_\si}\vv<\infty,$$ then by Theorem \ref{T4.1}, for any $\eta\in \C_\tau^1\cap \big(\cup_{n\ge 1}C([-\tau,0];\H_{A,n})\big)$ the form
$(\EE_\eta, C_b^2(\C_\tau))$ is closable on $L^2(\mu)$.

\end{document}